\documentclass[11pt]{article}
\usepackage{amsmath}
\usepackage{amssymb}
\usepackage{amsthm}
\usepackage{url}
\newcommand\ma{\medskipamount}
\newcommand\ba{\bigskipamount}
\newcommand\mLP{\\[\ma]}
\newcommand\mPP{\\[\ma]\indent}
\newcommand\bLP{\\[\ba]}
\newcommand\bPP{\\[\ba]\indent}
\newcommand\CC{\mathbb{C}}
\newcommand\RR{\mathbb{R}}
\newcommand\ZZ{\mathbb{Z}}
\newcommand\tha\theta
\newcommand\la\lambda
\newcommand\si\sigma
\newcommand\iy\infty
\newcommand\half{\frac12}
\newcommand\thalf{\tfrac12}
\newcommand\LHS{left-hand side}
\newcommand\RHS{right-hand side}
\renewcommand\Im{{\rm Im}\,}
\numberwithin{equation}{section}
\newtheorem{theorem}{Theorem}[section]
\newtheorem{Remark}[theorem]{Remark}
\newenvironment{remark}{\begin{Remark}\rm}{\end{Remark}}
\newcommand\Proof{\noindent{\bf Proof}\quad}

\begin{document}
\title{On the equivalence of two fundamental theta identities}
\author{Tom H. Koornwinder}
\date{\em Dedicated to the memory of Frank W. J. Olver}
\maketitle
%
%until 45
%
\begin{abstract}
Two fundamental theta identities, a three-term identity due to Weierstrass and a 
five-term identity due to Jacobi, both with products of four
theta functions as terms, are shown to be equivalent. One half of the equivalence was already
proved by R. J. Chapman in 1996. The history and usage of the two identities, and some
generalizations are also discussed.
\end{abstract}
\section{Introduction}
\label{36}
Theta functions occur in many parts of mathematics and its applications
\cite{25}.
While they had roots in the work of Jakob Bernoulli and Euler, they were introduced
in full
generality, depending on two arguments, by Jacobi. They became very
important in
nineteenth century complex analysis \cite{36}, \cite[Ch.~11]{27}
because elliptic functions could be expressed in terms of them.
Theta functions in several variables, later called {\em Riemann
theta functions} \cite[\S21.2]{3}, played a similar role for
abelian functions. Riemann's geometric approach \cite{37}
and Weierstrass' analytic approach \cite{38} were opposed to each other.
Ramanujan, in Chapters 16--21 of his Notebooks, gave many new theta identities,
with an emphasis on number theoretical and modular aspects; see Berndt \cite{47}
for an edited version with proofs.
Algebraic geometry, number theory and combinatorics are
some of the fields
where theta functions have played an important role since long.
New fields of application
arose during the last decades of the twentieth century:
nonlinear pde's like KdV \cite{28}, 
solvable models in statistical mechanics \cite{6},
Sklyanin algebra \cite{20}, \cite{21}, elliptic quantum groups \cite{34}
and elliptic hypergeometric series \cite{5}, \cite[Ch.~11]{2}, \cite{7}.

In literature identities involving theta functions abound, see for instance
Whittaker \& Watson \cite[Ch.~21]{1}, Erd\'elyi et al.\ \cite[\S13.10]{4} and
Olver et al.\ \cite[Ch.~20]{3},
but two identities (given in \cite{1} and \cite{3})
stand out because of their fundamental nature and because many of the other
identities can be derived from them. Both have the form of a sum of products
of four theta functions of different arguments being zero, with three terms in the first
formula and five terms in the second formula.
\mLP
{\bf First fundamental theta identity}
\begin{multline}
\tha_1(u+u_1)\tha_1(u-u_1)\tha_1(u_2+u_3)\tha_1(u_2-u_3)
+\tha_1(u+u_2)\tha_1(u-u_2)\tha_1(u_3+u_1)\tha_1(u_3-u_1)\\
+\tha_1(u+u_3)\tha_1(u-u_3)\tha_1(u_1+u_2)\tha_1(u_1-u_2)=0
\label{32}
\end{multline}
(or equivalently with $\tha_1$ replaced by $\si$), see
p.451, Example 5 and p.473, \S21.43 in \cite{1}, or
(23.10.4) and (23.6.9) in \cite{3}.
\mLP
{\bf Second fundamental theta identity}
\begin{multline}
2\,\tha_1(w)\,\tha_1(x)\,\tha_1(y)\,\tha_1(z)
=\tha_1(w')\,\tha_1(x')\,\tha_1(y')\,\tha_1(z')
+\tha_2(w')\,\tha_2(x')\,\tha_2(y')\,\tha_2(z')\\
-\tha_3(w')\,\tha_3(x')\,\tha_3(y')\,\tha_3(z')
+\tha_4(w')\,\tha_4(x')\,\tha_4(y')\,\tha_4(z'),
\label{33}
\end{multline}
where
\begin{equation}
\begin{split}
2w'=-w+x+y+z,&\qquad
2x'=w-x+y+z,\\
2y'=\;\;\,w+x-y+z,&\qquad
2z'=w+x+y-z
\end{split}
\label{37}
\end{equation}
and similar equivalent identities starting with $\tha_2$, $\tha_3$ or $\tha_4$ on the \LHS,
see \S21.22 in \cite{1}, or (21.6.6) and (21.2.9) in \cite{3}
\mPP
Identity \eqref{33} (the oldest one) was first given by
Jacobi \cite[p.507, formula (A)]{8};
this paper is based on notes made by Borchardt
of a course of Jacobi which were later annotated by Jacobi.
It first entered in Jacobi's lectures of 1835--1836 and he was so excited by the result that he completely changed his approach to elliptic
functions, using \eqref{33} as a starting point \cite[p.220]{36}.

Identity \eqref{32} was first obtained by Weierstrass \cite[(1.)]{23}.
For the proof he refers to
Schwarz \cite[Art.~38, formula (1.)]{10} (these are edited notes of lectures by Weierstrass).
Weierstrass \cite{23} mentions  that he first gave this formula in his lectures in 1862.
He emphasizes that \eqref{32} is essentially different from Jacobi's formulas \eqref{33}
and variants.

Some papers in the last decades have attributed these formulas to Riemann,
although without reference.
Frenkel \& Turaev \cite[pp.~171--172]{5} call formula \eqref{32}
{\em Riemann's theta identity}. Some later authors \cite[(3.4)]{24}, \cite[(5.3)]{13},
\cite[(6)]{7} also use this
terminology or speak about {\em Riemann's addition formula}.
As for \eqref{33}, Mumford \cite{35}, \cite[p.16]{12} calls it
{\em Riemann's theta relation}.
However, I have not been able to find formula \eqref{32} or \eqref{33}
in \cite{37} or elsewhere in Riemann's publications \cite{9}.

Formula \eqref{33} has a generalization \cite[Ch.~2, \S6]{12}, \cite[\S21.6(i)]{3}
to theta functions in several variables, which is called a
{\em generalized Riemann theta identity} by Mumford.
Weierstrass \cite{23}, \cite{39} gave a generalization of both \eqref{32} and
\eqref{33}, respectively, to the several variables case. It is not immediately clear how
the results in \cite{12} and \cite{39} are related.

The main purpose of this paper is to show
in Section \ref{34}
that \eqref{32} and \eqref{33} easily follow from
each other, and therefore can be considered to be equivalent identities.
We will work in the notation \cite[(11.2.1)]{2} for theta functions which is now common in
work on elliptic hypergeometric series. Its big advantage is that we have only one 
theta function instead of four different ones, by which lists of formulas can be greatly shrinked.
Another feature of this notation is that we work multiplicatively instead of additively.
Instead of double \mbox{(quasi-)periodicity} we have quasi-invariance under multiplication of the
independent variable by $q$. This notation is introduced in Section \ref{35}.
Some variants and applications of the two fundamental formulas are given in
Section \ref{38}. For completeness the elegant proofs by complex analysis of
the two fundamental formulas are recalled in Section \ref{39} and some other proofs are
mentioned. An Appendix compares some four-term theta identities
(mostly extending \eqref{32}).

During the revision of an earlier version of this paper my attention was called to work by
Wenchang Chu. Apparently unaware of the earlier occurrence of \eqref{32} in the
literature he posed it as a problem in the Monthly
to prove this formula (as a solution to a functional equation).
Several authors gave solutions \cite{46}. One of them, R.~J. Chapman, derived \eqref{32}
from \eqref{33}, essentially the same proof as my proof in Section \ref{34}.
Thus the present paper has become even more a survey paper than originally intended.
\bLP
{\bf Acknowledgements}\quad
I thank an anonymous referee for corrections of minor errors, for important additional material,
and for posing a question which brought me to writing the Appendix.
I also thank Bruce Berndt, Hjalmar Rosengren and Michael Schlosser
for providing additional references,
and Michael Somos for interesting comments.
\section{Preliminaries}
\label{35}
Let $q$ and $\tau({\rm mod}\;2\ZZ)$ be related by
$q=e^{i\pi\tau}$ and assume that $0<|q|<1$, or equivalently $\Im\tau>0$.
We will define and notate the {\em theta function} of {\em nome} $q$
as in Gasper \& Rahman \cite[(11.2.1)]{2}:
\begin{align}
\tha(w)=\tha(w;q)&:=(w,q/w;q)_\iy=
\prod_{j=0}^\iy(1-q^jw)(1-q^{j+1}/w)\qquad(w\ne0),
\label{1}\\
\tha(w_1,\ldots,w_k)&:=\tha(w_1)\ldots\tha(w_k).
\end{align}
By Jacobi's triple product identity \cite[(1.6.1)]{2} we have
\begin{equation}
\tha(w;q)=\frac1{(q;q)_\iy}\,\sum_{k=-\iy}^\iy (-1)^k q^{\half k(k-1)}w^k.
\label{22}
\end{equation}
Clearly,
\begin{align}
\tha(w^{-1};q)&=-w^{-1}\tha(w;q),\label{12}\\
\tha(qw;q)&=-w^{-1}\tha(w;q),\label{14}\\
\tha(q^k w;q)&=(-1)^k q^{-\half k(k-1)} w^{-k}\,\tha(w;q)\qquad(k\in\ZZ).
\label{13}
\end{align}

The four {\em Jacobi theta functions} $\tha_a$ or $\vartheta_a$
($a=1,2,3,4$), written as
\begin{equation*}
\tha_a(z)=\tha_a(z,q)=\tha_a(z\mid\tau)=\vartheta_a(\pi z,q)
=\vartheta_a(\pi z\mid\tau),
\end{equation*}
can all be expressed in terms of the theta function \eqref{1}:
\begin{align*}
\tha_1(z):=&i\,q^{1/4}(q^2;q^2)_\iy\,e^{-\pi iz}\,\tha(e^{2\pi iz};q^2),\\
\tha_2(z):=&q^{1/4}(q^2;q^2)_\iy\,e^{-\pi iz}\,\tha(-e^{2\pi iz};q^2)
=\tha_1(z+\thalf),\\
\tha_3(z):=&(q^2;q^2)_\iy\,\tha(-q\,e^{2\pi iz};q^2)
=\sum_{k=-\iy}^\iy q^{k^2}e^{2\pi ikz},\\
\tha_4(z):=&(q^2;q^2)_\iy\,\tha(q\,e^{2\pi iz};q^2)=\tha_3(z+\thalf).
\end{align*}
Note that
$\tha_1(z)$ is odd in $z$, while $\tha_2(z)$, $\tha_3(z)$ and
$\tha_4(z)$ are even in $z$.

The notation $\tha_a$ is used in \cite[\S13.10]{4} and
\cite[Ch.~20]{3}, while the notation $\vartheta_a$ is used in
\cite[Ch.~21]{1}. Mumford \cite{12} writes $\vartheta(z,\tau)$ instead of
$\theta_3(z\mid\tau)$.

The first fundamental identity \eqref{32} now takes the form
\begin{equation}
yu\,\tha(xy,x/y,vu,v/u)+uv\,\tha(xu,x/u,yv,y/v)+vy\,\tha(xv,x/v,uy,u/y)=0,
\label{2}
\end{equation}
or variants by applying \eqref{12}, see  \cite[(11.4.3)]{2}.
The terms in \eqref{2} are
obtained from each other by cyclic permutation in $y,u,v$.

The second fundamental identity \eqref{33} can be rewritten in the notation \eqref{1} as
\begin{multline}
2\tha(w^2,x^2,y^2,z^2;q^2)=\tha(w'',x'',y'',z'';q^2)+\tha(-w'',-x'',-y'',-z'';q^2)\\
+q^{-1}xyzw\big(\tha(qw'',qx'',qy'',qz'';q^2)-\tha(-qw'',-qx'',-qy'',-qz'';q^2)\big),
\label{16}
\end{multline}
where
\begin{equation*}
w''=w^{-1}xyz,\quad
x''=wx^{-1}yz,\quad
y''=wxy^{-1}z,\quad
z''=wxyz^{-1}.
\end{equation*}
\section{Variants and applications of the two fundamental formulas}
\label{38}
As already observed in Section \ref{36}, Weierstrass wrote \eqref{32} as
\begin{multline}
\si(u+u_1)\si(u-u_1)\si(u_2+u_3)\si(u_2-u_3)
+\si(u+u_2)\si(u-u_2)\si(u_3+u_1)\si(u_3-u_1)\\
+\si(u+u_3)\si(u-u_3)\si(u_1+u_2)\si(u_1-u_2)=0.
\label{28}
\end{multline}
The two formulas \eqref{32} and \eqref{28} are equivalent because
by \cite[p.473, \S21.43]{1}, for periods 1 and~$\tau$, we have
$\si(z)=C\,e^{\eta_1 z^2/2}\,\tha_1(\thalf z\mid\tau)$
with $C$ and $\eta_1$ only depending on~$\tau$.
For $v=1$, $u=-1$ formula \eqref{2} yields (using \eqref{12}):
\begin{equation}
\frac{y\,\tha(xy,x/y)}{\tha(x)^2\,\tha(y)^2}=f(y)-f(x),\quad\mbox{where }
f(x):=\frac{\tha(-x)^2}{\tha(-1)^2\,\tha(x)^2}\,.
\label{30}
\end{equation}
Conversely, \eqref{30} implies \eqref{2}, for any choice of the function $f$.
Indeed, after two substitutions from \eqref{30}, the first term of \eqref{2}
becomes
$\big(\tha(x)\tha(y)\tha(u)\tha(v)\big)^2\,(f(y)-f(x))(f(u)-f(v))$,
which adds up to 0 under cyclic permutation of $y,u,v$.

Weierstrass, according to Schwarz \cite[Art.~38, formula (1.)]{10}, derived \eqref{28} from the formula
\cite[Art.~11, formula (1.)]{10}:
\begin{equation}
\frac{\si(u+v)\si(u-v)}{\si^2(u)\si^2(v)}=\wp(v)-\wp(u),
\label{31}
\end{equation}
where $\wp(z)$ is Weierstrass' elliptic function.
By the expression \cite[13.20(4)]{4} of $\wp(z)$ in terms of theta functions,
\eqref{31} is equivalent to \eqref{30} just as \eqref{28} is equivalent to \eqref{2}.
Whittaker \& Watson give these results in
\cite[p.451, Examples 1 and 5]{1}. Formula \eqref{31} is also given in
\cite[(10.13.17)]{4} and \cite[(23.10.3)]{3}.

All addition formulas for theta functions in
\cite[pp.~487--488, Examples 1, 2, 3]{1} are instances of \eqref{30} or slight variants
of it which can be obtained by specialization of \eqref{2}.
Some of these formulas are used
in the proof that certain actions of the generators of the Sklyanin algebra on the space
of meromorphic functions determine a representation of the Sklyanin algebra
\cite[Theorem 2]{21}.

Weierstrass \cite{23}
observed at the end of his paper
that \eqref{28}, as a functional equation for
the sigma function, has a general solution given by a power series and
still depending on four arbitrary constants. This was finally proved in
full rigor by Hurwitz \cite{29}. However, \cite[pp. 452, 461]{1} gives earlier references for this
result to books by Halphen and by Hermite.

A stronger characterization of the sigma function (up to trivial transformations) and its degenerate cases was given by Bonk \cite{44}
as continuous solutions $\tau\colon\RR^n\to\CC$ of the functional equation
\[
\tau(u+v)\tau(u-v)=f_1(u)g_1(v)+f_2(u)g_2(v),
\]
where $f_1,f_2,g_1,g_2\colon\RR^n\to\CC$ are arbitrary.
Note that \eqref{31} is of the above form.

Elliptic, and in particular theta functions, entered in work
on solvable models in statistical mechanics started by
Baxter \cite{6} and followed up in papers like \cite{30}, \cite{31},
\cite{32}. While building on these publications, Frenkel \& Turaev
\cite{5} in their work on
the elliptic $6j$-symbol introduced elliptic hypergeometric series.
Among others, they obtained the summation formula of the terminating well-poised theta hypergeometric series ${}_{10}V_9(a;b,c,d,e,q^{-n};q,p)$. 
Formula \eqref{2} occurs as the first non-trivial case $n=1$ and it also plays a role
in the further proof by induction of this summation formula
\cite[\S11.4]{2}.
Closely related to these developments is the introduction of
elliptic quantum groups by Felder \cite{34}.
Again theta functions play here an important role \cite{33}, \cite{15}.
In \cite[Remarks 2.4, 4.3]{15} formula \eqref{2} is used in connection with the representation theory of the elliptic $U(2)$ quantum group.

If we pass in \eqref{2} to
homogeneous coordinates $(a_1,a_2,a_3,b_1,b_2,b_3)$ satisfying $a_1a_2a_3=b_1b_2b_3$
and expressed in terms of $x,y,u,v$ by
\[
a_1=b_3xu,\;
a_2=b_3xy,\;
a_3=b_3xv,\;
b_1=b_3x^2,\;
b_2=b_3xyuv, 
\]
then, after repeated application of \eqref{12}, we obtain another symmetric version
of \eqref{2}:
\begin{equation}
\frac{\tha(a_1/b_1,a_1/b_2,a_1/b_3)}{\tha(a_1/a_2,a_1/a_3)}
+\frac{\tha(a_2/b_1,a_2/b_2,a_2/b_3)}{\tha(a_2/a_3,a_2/a_1)}
+\frac{\tha(a_3/b_1,a_3/b_2,a_3/b_3)}{\tha(a_3/a_1,a_3/a_2)}=0\quad
(a_1a_2a_3=b_1b_2b_3).
\label{29}
\end{equation}
Formula \eqref{29}
has an $n$-term generalization which is associated with root system $A_{n-1}$:
\begin{equation}
\sum_{k=1}^n\frac{\prod_{j=1}^n \tha(a_k/b_j)}{\prod_{j\ne k}\tha(a_k/a_j)}=0\quad
(a_1\ldots a_n=b_1\ldots b_n),
\label{41}
\end{equation}
see \cite[(5)]{53},
\cite[Lemma A.2]{16}, \cite[(4.1)]{14}, \cite[Exercise 5.23]{2}.
The formula is given in terms of $\si(z)$ in
\cite[p.451, Example 3]{1}. Rosengren \cite[p.425]{14} traced the formula back to Tannery and Molk \cite[p.34]{17}.
Another $n$-term generalization, which reduces for $n=3$ to \eqref{2} after application of
\eqref{12}, is associated with root system $D_{n-1}$, see
\cite[Lemma 4.14]{18}, \cite[Lemma A.1]{16}, \cite[(4.6)]{14}.
More complicated many-term identities of theta products are given by
Kajihara \& Noumi \cite[Theorem 1.3]{54} and by
Langer, Schlosser \& Warnaar \cite[Theorem 1.1]{55}.

Of particular interest are  the four-term cases of the above identities.
Quite some four-term theta identities can be found on scattered places
in literature, and they also arise as special cases of
some identities in \cite{2}. One may wonder if some of these identities are essentially
the same (a question asked by a referee).
In one important case the answer will be negative, see the Appendix.

Various proofs and applications of \eqref{2} were given in \cite{46}, \cite{49}, \cite{50}.
Schlosser \cite[Remark 3.3]{52} observed that \eqref{2} pops up in connection
with a telescoping property of the special case $e=a^2/(bcd)$ of
Bailey's ${}_6\psi_6$ summation formula \cite[(II.33)]{2}.
This was converted into a proof of \eqref{2} by Chu \cite[Theorem 1.1]{48}.
In \cite[p.948]{24} formula \eqref{2} is used in the proof of a determinant evaluation
associated to the affine root system of type $C$.
In \cite{40} a $3\times 3$ determinant with theta function entries is evaluated, thus solving
an open problem in \cite{41}. The determinant evaluation has \eqref{2} as a special case.
In \cite{42} the so-called quintuple product identity
\cite[Exercise 5.6]{2} is derived from \eqref{2}.
\mPP
The second fundamental formula \eqref{33} and its variants can be written in a very
compact form by using the notation (cf.~\eqref{37})
\[
[a]:=\tha_a(w)\tha_a(x)\tha_a(y)\tha_a(z),\quad
[a]':=\tha_a(w')\tha_a(x')\tha_a(y')\tha_a(z').
\]
Then (the first one implies the others):
\begin{equation}
\begin{split}
2\,[1]=\;\;\,\,[1]'+[2]'-[3]'+[4]',&\quad
2\,[2]=[1]'+[2]'+[3]'-[4]',\\
2\,[3]=-[1]'+[2]'+[3]'+[4]',&\quad
2\,[4]=[1]'-[2]'+[3]'+[4]'.
\end{split}
\label{4}
\end{equation}
These are easily seen to be equivalent with
\cite[p.468, Example 1 and p.488, Example 7]{1}:
\begin{equation}
\begin{split}
[1]+[2]&=[1]'+[2]',\quad[1]+[3]=[2]'+[4]',\quad[1]+[4]=[1]'+[4]',\\
[1]-[2]&=[4]'-[3]',\quad[1]-[3]=[1]'-[3]',\quad[1]-[4]=[2]'-[3]'.
\end{split}
\label{3}
\end{equation}
Jacobi \cite[p.507, formula (A)]{8}
first obtained \eqref{3} and then derived \eqref{4} from it.

For $x=y=z=w$ \eqref{33} implies \cite[p.469, Example 4]{1}
\[
\tha_1(z)^4+\tha_3(z)^4=\tha_2(z)^4+\tha_4(z)^4.
\]

The computation \cite[Proposition 3]{21}
of the action of the Casimir operators in the representation of the Sklyanin algebra uses
\eqref{4}.
\section{Proofs of the fundamental theta relations}
\label{39}
For completeness I recall here the short and elegant complex analysis proofs
of the fundamental theta relations \eqref{2} and \eqref{16}.
\mLP
{\bf Proof of \eqref{2}} (Baxter \cite[p.460]{6}, see also \cite[p.3]{7}).\\
Consider the theta functions in \eqref{2} with nome $q^2$.
For fixed $y,u,v$ we have to prove that
\[
F(x):=\frac{yv^{-1}\,\tha(xy,x/y,vu,v/u;q^2)+yu^{-1}\,\tha(xv,x/v,uy,u/y;q^2)}
{\tha(xu,x/u,yv,y/v;q^2)}
\]
is equal to $-1$.
For generic values of $y,u,v$ $F(x)$ is a meromorphic function of $x$ on
$\CC\backslash\{0\}$. Then the numerator
vanishes at all (generically simple) zeros $x=q^{2k} u^{\pm1}$ ($k\in\ZZ$)
of the denominator. Indeed, for these values of $x$  the numerator equals
\begin{align*}
&yv^{-1}\,\tha(q^{2k} u^{\pm1}y,q^{2k} u^{\pm1}y^{-1},vu,vu^{-1};q^2)
+yu^{-1}\,\tha(q^{2k}u^{\pm1}v,q^{2k} u^{\pm1}v^{-1},uy,uy^{-1};q^2)\\
&=q^{-2k(k-1)}u^{\mp2k}\big(yv^{-1} \tha(u^{\pm1}y,u^{\pm1}y^{-1},vu,vu^{-1};q^2)
+yu^{-1}\,\tha(u^{\pm1}v,u^{\pm1}v^{-1},uy,uy^{-1};q^2)\big)=0,
\end{align*}
where we used \eqref{13} and \eqref{12}.
Thus $F$ is analytic in $x$ on $\CC\backslash\{0\}$.
Furthermore, $F(q^2x)=F(x)$ by \eqref{14}.
Hence $F$ is bounded.
Thus the singularity of $F$ at $0$ is removable and,
by Liouville's theorem, $F$ is constant. Now check that $F(v)=-1$
by \eqref{12}.\qed
\bPP
Whittaker \& Watson \cite[p.451, Examples 1 and 5]{1} obtain \eqref{2} from
\eqref{31}. They suggest a proof of \eqref{31}
by comparing zeros and poles of elliptic functions on both sides.
Liu \cite[(3.34)]{22} proves \eqref{2} by using a kind of generalized addition formula
for $\tha_1$.

Bailey \cite[(5.2)]{19} gives a more computational proof of \eqref{2}.
Among others he derives a three-term identity
\cite[(4.6)]{19} for very well-poised ${}_8\phi_7$ series, which Gasper \& Rahman
\cite[Exercise 2.15]{2} write in elegant symmetric form. 
By \cite[Exercise 2.16]{2}  formula \eqref{2} then should follow
from this three-term identity. Indeed, reduce it to a three-term identity
of very well-poised ${}_6\phi_5$ series which are summable by
\cite[(2.7.1)]{2}. See also \cite[Exercise 5.21]{2}.

Schlosser \cite[after (4.2)]{43} points out that \eqref{2} is also the special case
$b=1$ of \cite[(2.11.7)]{2} (put
$a=qxv$,  $c=xy$, $d=qx/y$, $e=qvu$, $f=qv/u$ and use \eqref{12} and \eqref{14}).
\bLP          
{\bf Proof of \eqref{16}} (Whittaker \& Watson \cite[p.468]{1}).\\
Divide the \RHS\ by the \LHS\ and consider the
resulting expression as a meromorphic function $F(w)$ of $w$ on
$\CC\backslash\{0\}$ (the other variables generically fixed)
with possible simple poles at the zeros $\pm q^k$ ($k\in\ZZ$) of
$\tha(w^2;q^2)$. Since $F(w)=F(-w)$ we can write $F(w)=G(w^2)$, where
$G$ is a meromorphic function on $\CC\backslash\{0\}$ with possible
simple poles at $q^{2k}$ ($k\in\ZZ$).
We have $F(qw)=F(w)$ because, by \eqref{14},
\begin{align*}
\frac{\tha(\pm q^{-1}w'';q^2)}{\tha(q^2w^2;q^2)}&=
\frac{\pm q^{-1}xyzw\,\tha(\pm qw'';q^2)}{\tha(w^2;q^2)}\,,\\
\frac{\pm xyzw\,\tha(\pm q^2x'',\pm q^2y'',\pm q^2z'';q^2)}{\tha(q^2w^2;q^2)}&=
\frac{\tha(\pm x'',\pm y'',\pm z'';q^2)}{\tha(w^2;q^2)}\,.
\end{align*}
Hence $G(q^2u)=G(u)$. But then
\[
2\pi i\,{\rm Res}_{u=q^{2k}} \big(u^{-1} G(u)\big)
=\int_{|u|=|q|^{2k-1}} G(u)\,\frac{du}u-
\int_{|u|=|q|^{2k+1}} G(u)\,\frac{du}u=0.
\]
Hence $G$ has no poles and similarly for $F$. Similarly as in the
previous proof we conclude that $F$ is constant in $w$. By symmetry,
$F$ is also constant in $x$, $y$ and $z$.
Thus we have shown that
\begin{multline}
A\,\tha(w^2,x^2,y^2,z^2;q^2)=\tha(w'',x'',y'',z'';q^2)+\tha(-w'',-x'',-y'',-z'';q^2)\\
+q^{-1}xyzw\big(\tha(qw'',qx'',qy'',qz'';q^2)-\tha(-qw'',-qx'',-qy'',-qz'';q^2)\big),
\label{18}
\end{multline}
for some constant $A$. Put in \eqref{18}
$w=x=q^\half$ and $y=z=iq$. Then $w''=x''=-q^2$ and $y''=z''=q$ and
\[
A\,\tha(q,q,-q^2,-q^2;q^2)=
\tha(q,q,-q^2,-q^2;q^2)+q^2\,\tha(q^3,q^3,-q^2,-q^2;q^2).
\]
Hence $A=2$ by \eqref{14}.\qed
\bPP
The last part of this proof is a slight improvement compared to
\cite[p.468]{1}. There it is first proved in
\cite[\S21.2]{1} (again by the same method) that
\begin{equation}
\tha(q;q^2)^2\,\tha(qz;q^2)^2=\tha(-q;q^2)^2\,\tha(-qz;q^2)^2
-qz\,\tha(-q^2;q^2)^2\,\tha(-q^2z;q^2)^2,
\label{20}
\end{equation}
and hence, by putting $z=1$,
\begin{equation}
\tha(q;q^2)^4=\tha(-q;q^2)^4-q\,\tha(-q^2;q^2)^4.
\label{21}
\end{equation}
Then the value of $A$ in the above proof is obtained by putting
$w=x=y=z=q^\half$ in \eqref{18} and comparing with \eqref{21}.

Note that \eqref{20} and \eqref{21} are special cases of \eqref{16}.
Another special case of \eqref{16}, only leaving two terms nonzero, was suggested by a referee:
\[
(w,x,y,z):=(i/z,q^{-1/2}/z,-iq^{1/2}/z,z).\quad\mbox{Then}\quad
(w'',x'',y'',z'')=(-1,q,-q^{-1},z^{-4}).
\]
Then \eqref{16} degenerates to
\begin{equation}
2\tha(-z^{-2},q^{-1}z^{-2},-qz^{-2},z^2;q^2)=\tha(-1,q,-q^{-1},z^{-4};q^2).
\label{40}
\end{equation}
This can be independently proved by \eqref{1} together with repeated use of the
standard identities
\[
(a,qa;q^2)_\iy=(a;q)_\iy\quad\mbox{and}\quad(a,-a;q)_\iy=(a^2;q^2)_\iy.
\]
Of course, \eqref{40} can also be used to settle that $A=2$ in \eqref{18}.

In Jacobi \cite[pp.~505--507]{8} and in Mumford \cite[Ch.~1, \S5]{12} a different proof
of \eqref{16} is given. It uses \eqref{22}.

If we compare our proofs of \eqref{2} and \eqref{16} given above with each other then we
see that in the proof of \eqref{2} it is not automatic that the possible simple poles have
residue zero because there are two simple poles in each annulus to be considered.
So we have to check there by computation that the numerator of $F(z)$ vanishes
whenever the denominator vanishes.
\section{Equivalence of the two fundamental theta relations}
\label{34}
Let us rewrite the first fundamental theta relation \eqref{2} as
$F_1(x,y,u,v;q)=0$, where
\begin{equation}
F_1(x,y,u,v;q):=
\tha(xy,x/y,uv,u/v;q^2)-\tha(xv,x/v,uy,u/y;q^2)-
uy^{-1}\tha(yv,y/v,xu,x/u;q^2).
\label{24}
\end{equation}
In the second fundamental theta relation \eqref{16} both sides are invariant under each
of the transformations of variable $w\to -w$, $x\to -x$, $y\to -y$, $z\to -z$.
Therefore we obtain an equivalent identity if we replace in \eqref{16}
$(w^2,x^2,y^2,z^2)$ by $(xy,x/y,uv,u/v)$. Thus we can write \eqref{16} equivalently, in
a form closer to \eqref{24}, as $F_2(x,y,u,v;q)=0$, where
\begin{multline}
F_2(x,y,u,v;q)\\:=2\tha(xy,x/y,uv,u/v;q^2)-\tha(xv,x/v,uy,u/y;q^2)
-\tha(-xv,-x/v,-uy,-u/y;q^2)\\
-q^{-1}xu\big(\tha(qxv,qx/v,quy,qu/y;q^2)-\tha(-qxv,-qx/v,-quy,-qu/y;q^2)\big).
\label{5}
\end{multline}

\begin{theorem}
The formulas $F_1(x,y,u,v;q)=0$ and $F_2(x,y,u,v;q)=0$ are equivalent to each other
because of the following identities:
\begin{align}
&F_1(x,y,u,v;q)+F_1(-x,y,-u,v;q)
-xy F_1(qx,qy,u,v;q)-xy F_1(-qx,qy,-u,v;q)\nonumber\\
&\qquad\qquad\qquad\qquad\qquad\qquad\qquad\qquad\qquad\qquad
\qquad\qquad\qquad\quad=F_2(x,y,u,v;q),\label{23}\\
&F_2(x,y,u,v;q)-uy^{-1} F_2(x,u,y,v;q)=2F_1(x,y,u,v;q).
\label{25}
\end{align}
\end{theorem}
\Proof
For the proof of \eqref{23} 
substitute \eqref{24} in the \LHS\ of \eqref{23}. Then
this \LHS\ becomes
\begin{align*}
&\tha(xy,x/y,uv,u/v;q^2)-xy\,\tha(q^2xy,x/y,uv,u/v;q^2)\\
&+\tha(-xy,-x/y,-uv,-u/v;q^2)-xy\,\tha(-q^2xy,-x/y,-uv,-u/v;q^2)\\
&-\tha(xv,x/v,uy,u/y;q^2)-\tha(-xv,-x/v,-uy,-u/y;q^2)\\
&+xy\big(\tha(qxv,qx/v,quy,q^{-1}u/y;q^2)+\tha(-qxv,-qx/v,-quy,-q^{-1}u/y;q^2)\big),
\end{align*}
which equals the \RHS\ of \eqref{23} because of \eqref{14} and \eqref{5}.

For the proof of \eqref{25}
substitute  \eqref{5} in the \LHS\ of \eqref{25}. Then
\begin{align*}
&2\tha(xy,x/y,uv,u/v;q^2)-2uy^{-1}\,\tha(xu,x/u,yv,y/v;q^2)\\
&-\tha(xv,x/v,uy,u/y;q^2)+uy^{-1}\,\tha(xv,x/v,uy,y/u;q^2)\\
&-\tha(-xv,-x/v,-uy,-u/y;q^2)+uy^{-1}\,\tha(-xv,-x/v,-uy,-y/u;q^2)\\
&-q^{-1}xu\big(\tha(qxv,qx/v,quy,qu/y;q^2)-\tha(qxv,qx/v,quy,qy/u;q^2)\big)\\
&+q^{-1}xu\big(\tha(-qxv,-qx/v,-quy,-qu/y;q^2)-\tha(-qxv,-qx/v,-quy,-qy/u;q^2)\big),
\end{align*}
which equals the \RHS\ of \eqref{25} because of \eqref{12}, \eqref{14} and \eqref{24}.\qed
\bPP
One half of the equivalence proof, i.e., essentially \eqref{25}, has been given before by
R.~J. Chapman \cite{46}.
\begin{remark}
It would be interesting to see if the above equivalence extends to theta functions in several
variables (cf.\  \cite{23}, \cite{39} and \cite[\S21.6(i)]{3}).
Similarly the question arises if for root systems $A_{n-1}$ and $D_{n-1}$
there is not only a first fundamental theta identity \cite{14} but also a second fundamental
identity, equivalent to the first one.
\end{remark}
\appendix
\section{Appendix: four-term theta identities}
Consider the case $n=4$ of \eqref{41}:
\begin{equation}
\sum_{k=1}^4\frac{\prod_{j=1}^4 \tha(a_k/b_j)}{\prod_{j\ne k}\tha(a_k/a_j)}=0\quad
(a_1\ldots a_4=b_1\ldots b_4).
\label{45}
\end{equation}
This can be seen as a four-term identity essentially depending on
six free variables with each term having seven theta factors in the
numerator and none in the denominator.
Slater \cite{53} mentions that \eqref{45} can be
rewritten as her formula (3), which is reproduced in
\cite[Exercise 5.22]{2}. As also mentioned in \cite{53}, replacement of
$b,g,h$ by $ab,ag,ah$, respectively, let $a$ disappear from this formula.
We are left with the four-term theta identity
\begin{multline}
b\,\tha(cb,db,eb,fb,g,h,g/h;q)
-b\,\tha(ch,dh,eh,fh,b,g,g/b;q)\\
=g\,\tha(cg,dg,eg,fg,b,h,b/h;q)
-h\,\tha(c,d,e,f,b/h,g/h,g/b;q)
\label{43}
\end{multline}
under the side condition
\begin{equation}
bcdefgh=q^2.
\label{44}
\end{equation}

A referee raised the question, to be answered negatively,
whether \eqref{45} is essentially the
same identity as the case $n=1$ of the elliptic ${}_{12}V_{11}$ 
Frenkel-Turaev \cite{5} identity, see \cite[(11.2.23)]{2}, or
equivalently (as pointed out by the same referee) as
the case $b=1$ of Bailey's four-term identity \cite[(2.12.9)]{2}
or \cite[(III.39)]{2} with side conditions there on p.58 or p.366,
respectively.

First consider the case $n=1$ of \cite[(11.2.23)]{2}.
While taking into account the side condition
$\la=qa^2/(bcd)$ as well as formulas (11.2.19) and (11.2.5) in \cite{2},
and after substitution of
$\tha(q;p)=-q\tha(q^{-1};p)$, the mentioned identity becomes the following
explicit four-term theta identity:
\begin{multline*}
1-\,\frac{\tha(b,c,d,e,f,(qa)^3/(bcdef);p)}
{\tha(qa/b,qa/c,qa/d,qa/e,qa/f,bcdef/(qa)^2;p)}
=\frac{\tha(qa,qa/(ef),(qa)^2/(bcde),(qa)^2/(bcdf);p)}
{\tha(qa/e,qa/f,(qa)^2/(bcdef),(qa)^2/(bcd);p)}\\
\times\left(1-\,\frac{\tha(qa/(cd),qa/(bd),qa/(bc),e,f,(qa)^3/(bcdef);p)}
{\tha(qa/b,qa/c,qa/d,(qa)^2/(bcde),(qa)^2/(bcdf),ef/(qa);p)}\right).
\end{multline*}
Replace $a$ by $q^{-1}a$. This will eliminate $q$. Next replace $p$ by $q$. We get:
\begin{multline*}
1-\,\frac{\tha(b,c,d,e,f,a^3/(bcdef);q)}
{\tha(a/b,a/c,a/d,a/e,a/f,bcdef/a^2;q)}
=\frac{\tha(a,a/(ef),a^2/(bcde),a^2/(bcdf);q)}
{\tha(a/e,a/f,a^2/(bcdef),a^2/(bcd);q)}\\
\times\left(1-\,\frac{\tha(a/(cd),a/(bd),a/(bc),e,f,a^3/(bcdef);q)}
{\tha(a/b,a/c,a/d,a^2/(bcde),a^2/(bcdf),ef/a;q)}\right).
\end{multline*}
Now multiply both sides by
$\tha(a/b,a/c,a/d,a/e,a/f,bcdef/a^2,a^2/(bcd);q)$ and simplify by using \eqref{12}.
We obtain:
\begin{multline}
\tha(a/b,a/c,a/d,a/e,a/f,bcdef/a^2,a^2/(bcd);q)
-\tha(b,c,d,e,f,a^3/(bcdef),a^2/(bcd);q)\\
=-a^{-2}bcdef\,\tha(a,a/(ef),a^2/(bcde),a^2/(bcdf),a/b,a/c,a/d;q)\\
-a^{-1}bcd\,\tha(a,a/(cd),a/(bd),a/(bc),e,f,a^3/(bcdef);q).
\label{42}
\end{multline}

Next consider \cite[(III.39)]{2} with side conditions given there on p.366.
Its case $b=1$ takes by use of \eqref{1} and \eqref{12} the form
\begin{multline*}
1-a^{-1}\,\frac{\tha(c,d,e,f,g,h;q)}{\tha(c/a,d/a,e/a,f/a,g/a,h/a;q)}
=-\la\,\frac{\tha(a^{-1},f/\la,g/\la,h/\la;q)}{\tha(\la,f/a,g/a,h/a;q)}\\
+\frac{\tha(a^{-1},f,g,h,\la c/a,\la d/a,\la e/a;q)}{\tha(\la,c/a,d/a,e/a,f/a,g/a,h/a;q)}\,,
\end{multline*}
or equivalently
\begin{multline*}
\tha(c/a,d/a,e/a,f/a,g/a,h/a,\la;q)-a^{-1}\,\tha(c,d,e,f,g,h,\la;q)\\
=-\la\,\tha(a^{-1},f/\la,g/\la,h/\la,c/a,d/a,e/a;q)
+\tha(a^{-1},f,g,h,\la c/a,\la d/a,\la e/a;q).
\end{multline*}
The side conditions take for $b=1$ the form
\begin{equation*}
\la=\frac{qa^2}{cde}\,,\qquad
h=\frac{a^3 q^2}{cdefg}\,.
\end{equation*}
By substitution of these in the identity, the variables $\la$ and $h$ disappear from the identity.
Now replace $e,b,g$ by $b,g,e$ (cyclic permutation) and use \eqref{12} and \eqref{13}.
We arrive at \eqref{42}.

Both \eqref{42} and \eqref{43} (with side condition) have six free
variables and seven theta factors in each term.
However, in \eqref{43} each term has three theta factors with the property
that the product of the arguments of two of them is equal to the argument
of the third. On the other hand, in none of the terms of \eqref{42} there
is such a simple relation between the arguments of the theta factors.
Of course, they must satisfy a relation, but this is a more complicated
polynomial relation involving the arguments of more than three of the
theta factors. Therefore, \eqref{42} and \eqref{43} cannot be matched with each other.

Let us finally consider some other four-term theta identities.
The case $n=4$ of the $D_{n-1}$ identity is a four-term identity depending on six free variables
with each term having ten theta factors.
A special case of both \cite[Theorem 1.3]{54} and \cite[Theorem 1.1]{55}
is the four-term theta identity \cite[(1.4)]{55} with four free variables and
each term having six theta factors.
As pointed out to me by Michael Schlosser, a further specialization of this identity to
two free variables yields identity (3.5) in Bouttier et al.\ \cite{45}
(take  $w=z$, $x=v=1/\alpha$ and $t=\alpha^3$ in \cite[(1.4)]{55}).
The identity in \cite{45} was
obtained in a statistical mechanical context and it was proved in a way quite similar
to Baxter's proof of \eqref{2} (see Section~\ref{39}).
As Hjalmar Rosengren pointed out to me, the $m=2$ cases of Pfaffian evaluations in
\cite[Remark 2.1]{51} yield four-term identities of theta products as well.

\quad\\
\begin{footnotesize}
\begin{quote}
{ T. H. Koornwinder, Korteweg-de Vries Institute, University of
 Amsterdam,\\
 P.O.\ Box 94248, 1090 GE Amsterdam, The Netherlands;

\vspace{\smallskipamount}
email: }{\tt T.H.Koornwinder@uva.nl}
\end{quote}
\end{footnotesize}

\end{document}